\documentclass[12pt, twoside, leqno]{article}

\usepackage{amsmath,amsthm}
\usepackage{amssymb}

\usepackage{enumerate}

\usepackage{graphicx}

\newtheorem{thm}{Theorem}[section]

\theoremstyle{definition}

\numberwithin{equation}{section}

\begin{document}


\baselineskip=17pt


\title{A theorem about vector fields with the ``proportional volume property"}

\author{Fabiano G. Brito\\
Centro de Matem\'atica, Computa\c{c}\~ao e Cogni\c{c}\~ao\\ 
Universidade Federal do ABC\\
09.210-170 Santo Andr\'e, Brazil\\
E-mail: fabiano.brito@ufabc.edu.br
\and 
Andr\'e O. Gomes\\
Departamento de Matem\'atica\\
Instituto de Matem\'atica e Estat\'{\i}stica da USP\\ 
Universidade de S\~ao Paulo\\
05508-090 S\~ao Paulo, Brazil\\
E-mail: gomes@ime.usp.br
\and
Robson M. Mesquita\\
Campus Universit\'ario de Arraias\\ 
Universidade Federal do Tocantins\\
77330-000 Tocantins, Brazil\\
E-mail: rbmat.ime@gmail.com}

\date{}

\maketitle


\renewcommand{\thefootnote}{}

\footnote{2010 \emph{Mathematics Subject Classification}: Primary .}

\footnote{\emph{Key words and phrases}: volume of vector fields}

\renewcommand{\thefootnote}{\arabic{footnote}}
\setcounter{footnote}{0}


\begin{abstract}
In this paper, we define a certain \textit{proportional volume property} for an unit vector field on a spherical domain in $\mathbb{S}^{3}$. We prove that the volume of these vector fields has an absolute minimum and this value is equal to the volume of the Hopf vector field. Some examples of such vector fields are given. We also study the minimum energy of solenoidal vector fields which coincides with a Hopf flow along the boundary of a spherical domain of $\mathbb{S}^{2k+1}$.
\end{abstract}

\section{Introduction} 
The \textit{volume} of an unit vector field $X$ on a compact and oriented Riemannian $n$-manifold $K$ can be defined as the volume of the submanifold $X(K)$ of the unit tangent bundle equipped with the restriction of the Sasaki metric. It is given by
\begin{align*}
	\mathrm{vol}(X)=\int_{K}\sqrt{\mathrm{det}(Id+(\nabla X)^{T}\nabla X)}d\nu
\end{align*}
where $d\nu$ is the volume element determined by the metric and $\nabla$ is the
Levi-Civita connection. On the other hand, the \textit{energy} of an unit vector field $X$ defined on $K$ is given by
\begin{align*}
	\mathcal{E}(X)=\frac{n}{2}\mathrm{vol}(K)+\frac{1}{2}\int_{K}\left\|\nabla X\right\|^{2}d\nu
\end{align*}
The following theorem is well known
\begin{thm}
$[\mathrm{GZ}]$ The unit vector fields of minimum volume on $\mathbb{S}^{3}$ are precisely the
Hopf vector fields, and no others.
\end{thm}
There is an analogue of Theorem 1.1 for the energy functional
\begin{thm} $[\mathrm{B}]$
Hopf vector fields on $\mathbb{S}^{3}$ minimize the functional $\mathcal{E}$ and
Hopf vector fields are the unique unit vector fields on $\mathbb{S}^{3}$ to minimize $\mathcal{E}$.
\end{thm}
The case where K is a submanifold with boundary of $\mathbb{S}^{2k+1}$ was treated first in [BGN] and the following general ``boundary version" was obtained
\begin{thm}
$[\mathrm{BGN}]$ Let $U$ be an open set of the (2k+1)-dimensional unit sphere $\mathbb{S}^{2k+1}$, let $K\subset U$ be a connected (2k+1)-submanifold with boundary of the sphere $\mathbb{S}^{2k+1}$ and let $\vec{v}$ be an unit vector field on $U$ which coincides with a Hopf flow H along
the boundary of $K$. Then 
\begin{eqnarray*}
	\mathrm{vol}(\vec{v})\!\!\!&\geq&\!\!\!\frac{4^{k}}{{2k \choose k}}\mathrm{vol}(K)
	\\
	\\
	\mathcal{E}(\vec{v})\!\!\!&\geq&\!\!\!(\frac{2k+1}{2}+\frac{k}{2k-1})\mathrm{vol}(K)
\end{eqnarray*}
\end{thm}
If $k=1$ we obtain that $\mathrm{vol}(\vec{v})\geq 2\mathrm{vol}(K)$, the volume of the Hopf vector field. In the proof of Theorem 1.2 was used the following important map $\varphi_{t}^{\vec{v}}:\mathbb{S}^{2k+1}\longrightarrow \mathbb{S}^{2k+1}(\sqrt{1+t^{2}})$ given by $\varphi_{t}^{\vec{v}}(x)=x+t\vec{v}(x)$ which was first used in the Milnor's paper [M]. 
\\
\\
In this work we continue to study the volume and energy of vector fields defined on submanifolds with boundary $K\subset \mathbb{S}^{3}$. In order to find similar results, we say that an unit vector field $\vec{v}$ satisfies the ``proportional volume property" if the following inequality hold for some $t>0$
\begin{align*}
\frac{\mathrm{vol}(\varphi_{t}^{\vec{v}}(K))}{\mathrm{vol}(\mathbb{S}^{3}(\sqrt{1+t^{2}}))}\geq\frac{\mathrm{vol}(K)}{\mathrm{vol}(\mathbb{S}^{3})}
\end{align*}
This condition is obviously equivalent to 
\begin{align*}
\mathrm{vol}(\varphi_{t}^{\vec{v}}(K))\geq\mathrm{vol}(K)(1+t^{2})^{3/2}
\end{align*}
We will show examples of vector fields satisfying the proportional volume property. Our goal is also to prove the following theorem

\begin{thm}
Let $K$ be a connected $3$-submanifold with boundary of the unit sphere $\mathbb{S}^{3}$ and let $\vec{v}$ be an unit vector field defined on $\mathbb{S}^{3}$. Consider the map $\varphi_{t}^{\vec{v}}:\mathbb{S}^{3}\longrightarrow \mathbb{S}^{3}(\sqrt{1+t^{2}})$ given by $\varphi_{t}^{\vec{v}}(x)=x+t\vec{v}(x)$ and a real number $t>0$ small enough so that the map $\varphi_{t}^{\vec{v}}$ is a diffeomorphism. Suppose also that the following conditions hold:
\begin{enumerate}
	\item $\displaystyle{\int_{\partial K}\left\langle \vec{v},\eta\right\rangle=0}$, where $\eta$ is the conormal vector to the boundary $\partial K$.
	\item $\displaystyle{\frac{\mathrm{vol}(\varphi_{t}^{\vec{v}}(K))}{\mathrm{vol}(\mathbb{S}^{3}(\sqrt{1+t^{2}}))}\geq\frac{\mathrm{vol}(K)}{\mathrm{vol}(\mathbb{S}^{3})}}$.
\end{enumerate}
Then $\mathrm{vol}(\vec{v}) \geq \mathrm{vol}(H)$ and $\mathcal{E}(\vec{v})\geq \mathcal{E}(H)$, where $H$ is the Hopf flow.
\end{thm}
Clearly, if the vector field $\vec{v}$ is solenoidal, the condition \textit{1} of Theorem 1.4 is satisfied. As a corollary of Proposition 1 in [BS], there is a lower limit for the energy of solenoidal fields defined on odd-dimensional Euclidean spheres

\begin{thm}
$[\mathrm{BS}]$ The Hopf vector fields has minimum energy among all
solenoidal unit vector fields on the sphere $\mathbb{S}^{2k+1}$.
\end{thm}

Using the same techniques of [BGN] we also prove the following boundary version of Theorem 1.5
\begin{thm}
Let $U$ be an open set of the $(2k+1)$-dimensional unit sphere $\mathbb{S}^{2k+1}$ and let $K\subset U$ be a connected $(2k+1)$-submanifold with boundary of the sphere $\mathbb{S}^{2k+1}$. Let $\vec{v}$ be a solenoidal unit vector field on $U$ which coincides with a Hopf flow $H$ along the boundary of K. Then 
\begin{eqnarray*}
	\mathcal{E}(\vec{v})\geq \left(\frac{2k+1}{2}+k\right)\mathrm{vol}(K)= \mathcal{E}(H)
\end{eqnarray*}
\end{thm}

\section{Preliminaries}
Using an orthonormal local frame $\left\{e_{1},\ldots, e_{n-1}, e_{n}=\vec{v}\right\}$ of $K$, the volume of the unit vector field $\vec{v}$ is given by

\begin{eqnarray*}
\mathrm{vol}(\vec{v})=\int_{K} (1+\sum\limits_{a=1}^{n}\left\|\nabla_{e_{a}}\vec{v}\right\|^{2}+\sum\limits_{a<b}\left\|\nabla_{e_{a}}\vec{v}\wedge\nabla_{e_{b}}\vec{v}\right\|^{2}+\ldots
\\
\ldots+\sum\limits_{a_{1}<\cdots<a_{n-1}}\left\|\nabla_{e_{a_{1}}}\vec{v}\wedge\cdots\wedge\nabla_{e_{a_{n-1}}}\vec{v}\right\|^{2})^{1/2}\: d\nu
\end{eqnarray*}

Now, let $U\subset\mathbb{S}^{2k+1}$ be an open set of the unit sphere and let $K\subset U$ be a connected $(2k+1)$-submanifold with boundary of $\mathbb{S}^{2k+1}$. Suppose that $\vec{v}$ is an unit vector field defined on $U$. We also consider the map $\varphi_{t}^{\vec{v}}:U\longrightarrow \mathbb{S}^{2k+1}(\sqrt{1+t^{2}})$ given by $\varphi_{t}^{\vec{v}}(x)=x+t\vec{v}(x)$. We assume that $t>0$ is small enough so that the map $\varphi_{t}^{\vec{v}}$ is a diffeomorphism. 
\\
\\
In [BGN], the Authors showed that the determinant of the Jacobian matrix of $\varphi_{t}^{\vec{v}}$ can be express in the form 
\begin{eqnarray*}
	\det(d\varphi_{t}^{\vec{v}})=\sqrt{1+t^{2}}(1+\sum\limits_{i=1}^{2k}\sigma_{i}(\vec{v})t^{i})
\end{eqnarray*}
where, by definition, $h_{ij}(\vec{v}):=\left\langle \nabla_{e_{i}}\vec{v},e_{j}\right\rangle$ (with $i,j\in\left\{1,\ldots,2k\right\}$) and the functions $\sigma_{i}$ are the $i$-symmetric functions of the $h_{ij}$. For instance, if $k=1$
\begin{eqnarray*}
\sigma_{1}(\vec{v})\!\!\!&=&\!\!\!h_{11}(\vec{v})+h_{22}(\vec{v})\\ \sigma_{2}(\vec{v})\!\!\!&=&\!\!\!h_{11}(\vec{v})h_{22}(\vec{v})-h_{12}(\vec{v})h_{21}(\vec{v})
\end{eqnarray*}
and the volume $\mathrm{vol}(\vec{v})$ has the following form
\begin{align}
	\mathrm{vol}(\vec{v})=\int_{K}\left(\sqrt{1+\sum\limits_{i,j=1}^{2}h_{ij}^{2}+(\mathrm{det}(h_{ij}))^{2}+\cdots}\right)\: d\nu
\end{align}
For an arbitrary integer $k \geq 1$ the energy $\mathcal{E}(\vec{v})$ has the following form 
\begin{align}
	\mathcal{E}(\vec{v})=\frac{2k+1}{2}\mathrm{vol}(K)+\frac{1}{2}\int_{K}[\sum\limits_{i,j=1}^{2k}(h_{ij}(\vec{v}))^{2}+\sum\limits_{i=1}^{2k}(\left\langle\nabla_{\vec{v}}\vec{v},e_{i}\right\rangle)^{2}]
\end{align}
and, by change of variables theorem, we have the following expression for the volume of the spherical domain $\varphi_{t}^{\vec{v}}(K)\subset \mathbb{S}^{3}(\sqrt{1+t^{2}})$ 
\begin{eqnarray}	
\mathrm{vol}[\varphi_{t}^{\vec{v}}(K)]=\int_{K}\sqrt{1+t^{2}}(1+\sum\limits_{i=1}^{2k}\sigma_{i}(\vec{v})t^{i})
\end{eqnarray}

\section{Proof of theorem 1.4} 

By hypothesis \textit{2}) and by equality (2.3), we have 
\begin{align}
\mathrm{vol}[\varphi_{t}^{\vec{v}}(K)]=\int_{K}\sqrt{1+t^{2}}(1+\sigma_{1}(\vec{v})t+\sigma_{2}(\vec{v})t^{2})\geq \int_{K}(\sqrt{1+t^{2}})^{3}
\end{align}
On the other hand, by hypothesis \textit{1}) and by divergence theorem we have $\displaystyle{\int_{K}\sigma_{1}(\vec{v})=0}$. Then (3.1) has the form
\begin{align}
\int_{K}(1+\sigma_{2}(\vec{v})t^{2})\geq \int_{K}(1+t^{2})
\end{align}
and then
\begin{align}
	\int_{K}\sigma_{2}(\vec{v})\geq\mathrm{vol}(K)
\end{align}
Now, observing that
\begin{align*}
	\mathrm{vol}(H)=2\mathrm{vol}(K) \ \ \ \ \& \ \ \ \ \sum\limits_{i,j=1}^{2}h_{ij}^{2}\geq 2\sigma_{2}(\vec{v})
\end{align*}
we have by equation (2.1) that
\begin{align*}
	\mathrm{vol}(\vec{v})&=\int_{K}\sqrt{1+\sum h_{ij}^{2}+(\mathrm{det}(h_{ij}))^{2}+\ldots}\\
	&\geq \int_{K}\sqrt{1+2\sigma_{2}(\vec{v})+\sigma_{2}(\vec{v})^{2}}\\
	&\geq \int_{K}(1+\sigma_{2}(\vec{v}))\geq 2\mathrm{vol}(K)=\mathrm{vol}(H)
\end{align*}
In a similar way, using the equation (2.2) for $k=1$ we can write
\begin{align*}
	\mathcal{E}(\vec{v})\geq \frac{3}{2}\mathrm{vol}(K)+\frac{1}{2}\int_{K}\sum\limits_{i,j=1}^{2}(h_{ij}(\vec{v}))^{2}
\end{align*}
and then 
\begin{align*}
	\mathcal{E}(\vec{v})\geq\frac{3}{2}\mathrm{vol}(K)+\int_{K}\!\!\sigma_{2}(\vec{v})\geq \frac{3}{2}\mathrm{vol}(K)+\mathrm{vol}(K)=\mathcal{E}(H)
\end{align*}
\section{Proof of theorem 1.6}
Once again, as a consequence of equation (2.2) we have
\begin{eqnarray}
\mathcal{E}(\vec{v})\geq \frac{2k+1}{2}\mathrm{vol}(K)+\frac{1}{2}\int_{K}\sum\limits_{i,j=1}^{2k}(h_{ij}(\vec{v}))^{2}
\end{eqnarray}
Now observe that
\begin{eqnarray}
\sum\limits_{i<j}(h_{ii}-h_{jj})^{2}=(2k-1)\sum\limits_{i}h_{ii}^{2}-2\sum\limits_{i<j}h_{ii}h_{jj}
\end{eqnarray}
and as $\vec{v}$ is a solenoidal vector field
\begin{eqnarray}
	0=[\sigma_{1}(\vec{v})]^{2}=(\sum\limits_{i}h_{ii})^{2}=\sum\limits_{i}h_{ii}^{2}+2\sum\limits_{i<j}h_{ii}h_{jj}
\end{eqnarray}
in other words
\begin{eqnarray}
	-2\sum\limits_{i<j}h_{ii}h_{jj}=\sum\limits_{i}h_{ii}^{2}
\end{eqnarray}
Substituting equation (4.4) in (4.2) we obtain
\begin{eqnarray}
\sum\limits_{i<j}(h_{ii}-h_{jj})^{2}=-4k\sum\limits_{i<j}h_{ii}h_{jj}
\end{eqnarray}
Further, we also have the following equation
\begin{eqnarray}
\sum\limits_{i<j}(h_{ij}+h_{ji})^{2}=\sum\limits_{i\neq j}h_{ij}^{2}+2\sum\limits_{i<j}h_{ij}h_{ji}
\end{eqnarray}
and then
\begin{eqnarray}
2k\sum\limits_{i<j}(h_{ij}+h_{ji})^{2}=2k\sum\limits_{i\neq j}h_{ij}^{2}+4k\sum\limits_{i<j}h_{ij}h_{ji}
\end{eqnarray}
Adding equations (4.5) and (4.7), we have
\begin{eqnarray}
\sum\limits_{i\neq j}h_{ij}^{2}\geq 2\sigma_{2}
\end{eqnarray}
and
\begin{eqnarray}
\sum\limits_{i,j=1}^{2k}h_{ij}^{2}=\sum\limits_{i}h_{ii}^{2}+\sum\limits_{i\neq j}h_{ij}^{2}\geq 2\sigma_{2}
\end{eqnarray}
Using the inequalities (4.1) and (4.9), we find 
\begin{eqnarray}
\mathcal{E}(\vec{v})&\geq&\!\!\! \frac{2k+1}{2}\mathrm{vol}(K)+\int_{K}\sigma_{2}(\vec{v})
\end{eqnarray}
On the other hand, by a similar argument used in [BGN] we can say that
\begin{eqnarray}
	\int_{K}\sigma_{2}(\vec{v})=k\mathrm{vol}(K)
\end{eqnarray}
and then we have
\begin{eqnarray*}
	\mathcal{E}(\vec{v})\geq \frac{2k+1}{2}\mathrm{vol}(K)+k\mathrm{vol}(K)=\left(\frac{2k+1}{2}+k\right)\mathrm{vol}(K)
\end{eqnarray*}

\section{Examples}

Let $K\subset \mathbb{S}^{3}$ be a solid torus whose boundary is the Clifford torus $T$. We can take the following parametrization for $K$ 
	\begin{align*}
		x(\theta, \alpha, \delta)=\left(\delta\cos(\theta),\delta\sin(\theta),\sqrt{1-\delta^{2}}\cos(\alpha),\sqrt{1-\delta^{2}}\sin(\alpha)\right)
	\end{align*}
with $0\leq \theta,\alpha \leq 2\pi$ and $0\leq \delta \leq \frac{1}{\sqrt{2}}$.

Consider an unit vector field $\vec{v}$ defined on $\mathbb{S}^{3}$ and tangent to the $\partial K$. For example, an element of the family of solenoidal unit vector fields
\begin{align*}
	\vec{v}_{\lambda}(x,y,z,w)=\frac{1}{\sqrt{1+(\lambda^{2} -1)(x^{2}+y^{2})}}\left(-\lambda y,\lambda x,-w,z\right)
\end{align*}
where $\lambda >1$. Then it satisfies the first condition of Theorem 1.4. For the second condition, we have two possibilities: 
\begin{align*}
 \frac{\mathrm{vol}(\varphi_{t}^{\vec{v}}(K))}{\mathrm{vol}(\mathbb{S}^{3}(\sqrt{1+t^{2}}))}\geq\frac{\mathrm{vol}(K)}{\mathrm{vol}(\mathbb{S}^{3})}
\end{align*}
or 
 \begin{align*}
	\frac{\mathrm{vol}(\varphi_{t}^{\vec{v}}(K))}{\mathrm{vol}(\mathbb{S}^{3}(\sqrt{1+t^{2}}))}<\frac{\mathrm{vol}(K)}{\mathrm{vol}(\mathbb{S}^{3})}
\end{align*}
In the first alternative, we find the example. In the second alternative, we consider the complementar of $K$, that is, the solid torus $K^{c}\subset\mathbb{S}^{3}$ with boundary $T$ such that $K \cup K^{c}=\mathbb{S}^{3}$ and $K\cap K^{c}=T$. 
\\
\\
We claim that: $\displaystyle{\frac{\mathrm{vol}(\varphi_{t}^{\vec{v}}(K^{c}))}{\mathrm{vol}(\mathbb{S}^{3}(\sqrt{1+t^{2}}))}\geq\frac{\mathrm{vol}(K^{c})}{\mathrm{vol}(\mathbb{S}^{3})}}$. 
\\
\\
In fact, if $\displaystyle{\frac{\mathrm{vol}(\varphi_{t}^{\vec{v}}(K))}{\mathrm{vol}(\mathbb{S}^{3}(\sqrt{1+t^{2}}))}<\frac{\mathrm{vol}(K)}{\mathrm{vol}(\mathbb{S}^{3})}}$ and $\displaystyle{\frac{\mathrm{vol}(\varphi_{t}^{\vec{v}}(K^{c}))}{\mathrm{vol}(\mathbb{S}^{3}(\sqrt{1+t^{2}}))}<\frac{\mathrm{vol}(K^{c})}{\mathrm{vol}(\mathbb{S}^{3})}}$
\\
\\
then, adding the two inequalities we arrive at a contradiction
\begin{align*}
1=\frac{\mathrm{vol}(\varphi_{t}^{\vec{v}}(K))+\mathrm{vol}(\varphi_{t}^{\vec{v}}(K^{c}))}{\mathrm{vol}(\mathbb{S}^{3}(\sqrt{1+t^{2}}))}<\frac{\mathrm{vol}(K)+\mathrm{vol}(K^{c})}{\mathrm{vol}(\mathbb{S}^{3})}=1
\end{align*}
Therefore we find an example after substituting $K$ for $K^{c}$.

\subsection*{Dedication}
Dedicated to Professor Antonio Gervasio Colares on his 80th birthday.

\end{document}